\documentclass[reqno,12pt, oneside]{amsart}

\usepackage{enumerate}
\usepackage[nobysame]{amsrefs}
\usepackage{enumitem}
\usepackage{amssymb}
\usepackage{amsmath}
\usepackage{amsthm}
\usepackage{microtype}
\usepackage[right=2cm,left=2cm, top=2.5cm, bottom=2.5cm,includehead]{geometry}
\usepackage[mathlines]{lineno}

\newcommand*\patchAmsMathEnvironmentForLineno[1]{%
  \expandafter\let\csname old#1\expandafter\endcsname\csname #1\endcsname
  \expandafter\let\csname oldend#1\expandafter\endcsname\csname end#1\endcsname
  \renewenvironment{#1}%
     {\linenomath\csname old#1\endcsname}%
     {\csname oldend#1\endcsname\endlinenomath}}% 
\newcommand*\patchBothAmsMathEnvironmentsForLineno[1]{%
  \patchAmsMathEnvironmentForLineno{#1}%
  \patchAmsMathEnvironmentForLineno{#1*}}%
\AtBeginDocument{%
\patchBothAmsMathEnvironmentsForLineno{equation}%
\patchBothAmsMathEnvironmentsForLineno{align}%
\patchBothAmsMathEnvironmentsForLineno{flalign}%
\patchBothAmsMathEnvironmentsForLineno{alignat}%
\patchBothAmsMathEnvironmentsForLineno{gather}%
\patchBothAmsMathEnvironmentsForLineno{multline}%
}

\theoremstyle{plain}
\newtheorem{theorem}{Theorem}
\newtheorem{corollary}{Corollary}
\newtheorem{lemma}{Lemma}
\newtheorem{proposition}{Proposition}

\newcommand{\bn}{\mathbb{N}}

\newcommand{\br}{\mathbb{R}}

\theoremstyle{remark}

\newtheorem{example}{Example}

\theoremstyle{definition}
\newtheorem{definition}{Definition}

\DeclareMathOperator{\var}{var}
\DeclareMathOperator{\ivar}{ivar}

\newcommand{\norm}[1]{\left\|#1\right\|}

\newcommand{\abs}[1]{\lvert#1\rvert}
\newcommand{\babs}[1]{\bigl|#1\bigr|}

\newcommand{\zdef}{{\mathrel{\mathop:}=}}

\newcommand{\lnorm}[1]{\left\|#1\right\|_{\Lambda BV}}

\newcommand{\LBV}{{\Lambda BV}}

\begin{document}

\title[Compactness]{Compactness in the spaces of functions of bounded variation}

\author{Jacek Gulgowski}

\address{Institute of Mathematics, Faculty of Mathematics, Physics and Informatics, University of Gda\'nsk, 80-308 Gda\'nsk, Poland }

\email[J.~Gulgowski]{dzak@mat.ug.edu.pl}


\begin{abstract}
Recently the characterization of the compactness in the space $BV([0,1])$ of  functions of bounded Jordan variation was given. Here, certain generalizations of this result are given for the spaces of functions of bounded Waterman $\Lambda$-variation, Young $\Phi$-variation and integral variation. It appears that on the compact sets  the norm is uniformly approximated by certain seminorms induced by the selection of finitely many intervals in $[0,1]$.
\end{abstract}

\maketitle

\section{Preliminaries}
\label{sec:prem}

\subsection{Notation and basic definitions}

By $\mathbb N$ we denote the set of all positive  integers and  by $I$ we denote the unit interval $[0,1]$. We will call the family $(I_n)$ (countable, finite or infinite) of closed intervals $I_n\subset I$ a family of {\em nonoverlapping intervals} when for any two intervals $I_n,I_m$ the intersection $I_n\cap I_m$ consists of at most one point (their interiors are pairwise disjoint).

As usual, for $q \in [1,+\infty)$, by $L^q(J)$ we will denote the Banach space of all the equivalence classes of real-valued functions defined on a bounded interval $J \subseteq \mathbb R$ which are  Lebesgue integrable with $q$-th power, endowed with the norm $\norm{x}_{L^q} \zdef \bigl(\int_J \abs{x(t)}^q\textup dt\bigr)^{1/q}$. 
Further on, we will simply write  $L^q$ instead of $L^q(I)$. From now on, we will always assume that $q'\in(1,+\infty)$ is a conjugate of $q\in(1,+\infty)$, that is $1/q'+1/q = 1$. Lebesgue measure of the set $E\subset I$ will be denoted as $|E|$.

Below, the norm in a normed space $E$ will be denoted by $\norm{\cdot}_E$. The open ball centered at $x$ with  radius $r>0$, in a space $E$, will be denoted by $B_E(x,t)$.

The concept of a variation of the function has evolved since the end of 19th century when the definition was  introduced by Camille Jordan. Among different generalizations we should mention 
\begin{itemize}
\item the Wiener $p$-variation, for $p\in[1,+\infty)$;
\item the Young $\phi$-variation (defined for an appropriate function $\phi\colon[0,+\infty)\to\br$);
\item the Waterman $\Lambda$-variation (for an appropriate sequence $\Lambda = (\lambda_n)_{n\in\bn}\subset (0,+\infty)$);
\item  $q$-integral $p$-variation (for $1\leq p<q\leq +\infty$, introduced in 1970s in papers of Brudnyi and Terehin).
\end{itemize} 

Now we are going to collect  basic definitions and facts concerning functions of bounded variation which will be needed in the sequel.

\begin{definition}
Let $x\colon I\to \br$ be a real-valued function defined on  $I$. The number
\[
 \var x=\sup \sum_{i=1}^{N}\abs{x(t_i)-x(t_{i-1})},
\]
where the supremum is taken over all finite partitions $0=t_0<t_1<\ldots<t_N=1$ of $I$, is called the \emph{Jordan variation} (or just a variation, for short) of the function $x$ over $I$. If the interval $I$ is replaced with some other interval $[a,b]$ the variation of the function $x\colon[a,b]\to\br$ will be denoted by $\var(x,[a,b])$.
\end{definition}

\subsection{Waterman $\Lambda$-variation}

Let us first define the $\Lambda$-variation of the function $x\colon I\to\br$. The concept was introduced by Waterman in \cite{Wat1}. Since then the functions of bounded $\Lambda$-variation were intensely studied by many authors -- for an overview we refer to \cite{ABM}.

 \begin{definition}\label{LambdaSequence}
  Let us consider a nondecreasing sequence of positive real numbers $\Lambda = (\lambda_n)_{n\in\mathbb{N}}$.
  We call such sequence a \emph{Waterman sequence} if
  \[
   \sum_{n=1}^\infty\frac{1}{\lambda_n} = +\infty.
  \]
 \end{definition}
 
 \begin{definition}
  Let $\Lambda = (\lambda_n)_{n\in\mathbb{N}}$ be a Waterman sequence and let $x\colon I\to\mathbb{R}$. We say that $x$ is of bounded $\Lambda$-variation  if there exists a positive  constant $M$ such that for any finite sequence of nonoverlapping  subintervals $\{[a_1, b_1], [a_2, b_2], \ldots, [a_N, b_N]\}$ of $I$, the following inequality holds
  \[
   \sum_{i=1}^N \frac{\abs{x(b_i) - x(a_i)}}{\lambda_i} \leq M.
  \]
  The supremum of the above sums, taken over the family of all sequences of nonoverlapping subintervals of $I$, is called the $\Lambda$-variation of $x$  and it is denoted by $\var_\Lambda(x)$. 
 \end{definition}
 
It is worth to mention that there are many equivalent ways to express that the function $x\colon I\to\br$ is of bounded $\Lambda$-variation (cf. \cite{Wat}*{Theorem 1, p.~34} and \cite[Lemma 1]{BCGS}), but we will not go into the details here.

 The space of all  functions defined on the interval $I$ and  of bounded $\Lambda$-variation, endowed with the norm $\lnorm{x} \zdef \abs{x(0)} + \var_\Lambda(x)$ forms a Banach space $\LBV(I)$ (see \cite{Wat}*{Section 3}).

We will also need a notation for finite sums approximating the $\Lambda$-variation. Let ${\mathcal I}=\{I_1, I_2, \ldots, I_N\}$ be a finite sequence of nonoverlapping   subintervals of $I$ in the form $I_i=[a_i,b_i]$  and let $x \colon I \to \mathbb R$ be a bounded function. Let us define
the  seminorm induced by ${\mathcal I}$:

\begin{equation}
\sigma_\Lambda(x,{\mathcal I}) \zdef  \sum_{i=1}^N \frac{\abs{x(b_i) - x(a_i)}}{\lambda_i}.
\end{equation}

\subsection{Young $\phi$-variation.}

Now we are going to remind definitions related to the Young $\phi$-variation.

\begin{definition}
A function $\phi \colon \mathbb [0,+\infty) \to \mathbb [0,+\infty)$ is said to be a \emph{$\phi$-function} if it is continuous, unbounded, nondecreasing, convex and such that $\phi(u)=0$ if and only if $u=0$.
\end{definition}

Using the notion of a $\phi$-function one can define the variation in the sense of Young.

\begin{definition}
Let $x$ be a real-valued function defined on $I$ and let $\phi$ be a given $\phi$-function. The number
\[
 \var_{\phi}x=\sup \sum_{i=1}^{N}\phi(\abs{x(t_{i})-x(t_{i-1})}),
\]
where the supremum is taken over all finite partitions $0=t_0<t_1<\ldots<t_N=1$ of $I$, is called the \emph{$\phi$-variation} (or \emph{variation in the sense of Young}) of the function $x$ over $I$.
\end{definition}

The set of all functions of bounded $\phi$-variation does not necessarily form a linear space, but it may naturally be extended to the 
 space $\Phi BV$ given by
\[
\Phi BV = \{ x : \exists_{\lambda >0}   \var_{\phi}(x/\lambda) < +\infty \}
\]
with the norm
\[
\Vert x\Vert_{\Phi BV} = |x(0)| + \inf\{ \lambda >0 : \var_{\phi}(x/\lambda) \leq 1\}.
\]
The space $\Phi BV$ with the norm $\Vert \cdot\Vert_{\Phi BV} $ forms a Banach space. Later we will denote the seminorm appearing in the formula above as
\[
V_\phi(x) = \inf\{ \lambda >0 : \var_{\phi}(x/\lambda) \leq 1\}.
\]

Similarly as before, also in case of $\Lambda$-variation, we may define the appropriate seminorms approximating the norm and based on the finite family  ${\mathcal I}=\{I_1, I_2, \ldots, I_N\}$ of nonoverlapping subintervals of $I$. We will need the following notation:

\begin{equation}
\sigma_\phi(x,{\mathcal I}) \zdef  \sum_{i=1}^N \phi(\abs{x(b_i) - x(a_i)});
\end{equation}

\begin{equation}
S_\phi(x,{\mathcal I}) \zdef  \inf\{ \lambda >0 : \sigma_{\phi}(x/\lambda,{\mathcal I}) \leq 1\}.
\end{equation} 
It is clear that for any family ${\mathcal I}$ of nonoverlapping intervals there is $S_\phi(x,{\mathcal I})\leq V_{\phi}(x)$.

\subsection{$q$-integral $p$-variation.}

Apart from the concepts given above we will look at the integral variation. The concept of the $q$-integral $p$-variation was introduced by Terehin (see \cite{T1}) and later studied by Borucka\--Cie\'{s}lewicz \cite{BC1,BC2}. Recently it was also investigated in \cite{JG:ibv}. We should  mention here that this definition was also investigated for multivariate maps from the approximation theory perspective: the concepts appeared first in \cite{yb71} and quite recently in \cite{Br} (here the partitions were made of cubes) or in \cite{kopotun} (here the partitions were made by {\em rings}, i.e. differences of two cubes). We should also point out that the  $q$-integral $p$-variation considered in this paper is a special case of the multivariate definition given in \cite[Definition 1.2]{BrBrDiss}.

Let us first remind the classical definition of the $L^q$ modulus of continuity (see e.g. \cite{Z}).  
\begin{definition}Let $x\colon I\to \br$ be a Lebesgue measurable function and let $[a,b]\subset I$ be a fixed interval, $a<b$. The value
\[
\omega_q(x; a,b) = \sup_{0 < h< b-a}\limits \Bigl( \int_a^{b-h} |x(t+h)-x(t)|^q\textup dt\Bigr)^{1/q} = \sup_{0 < h< b-a}\limits \Vert x(\cdot + h) - x(\cdot) \Vert_{L^q(a,b-h)}
\]
is called the $L^q$-modulus of continuity of a function $x$ on an interval $[a,b]$.  
\end{definition}

We may observe that if $x\in L^q(a,b)$ then its $L^q$-modulus of continuity is well defined and
\[
 \Vert x(\cdot + h) - x(\cdot) \Vert_{L^q(a,b-h)} \leq  \Vert x(\cdot + h) \Vert_{L^q(a,b-h)} +  \Vert x \Vert_{L^q(a,b-h)} \leq  2\Vert x\Vert_{L^q(a,b)}
\]
and
\begin{equation}\label{eq:Lqmod:est}
\omega_q(x; a,b) \leq 2\Vert x\Vert_{L^q(a,b)}.
\end{equation}

Let us now define the integral variation of a measurable function (cf. \cite{JG:ibv}).
\begin{definition}\label{def:ivar}
Let $p,q \in [1,+\infty)$ and let $x \colon [a,b] \to \mathbb R$ be a Lebesgue measurable function, where $[a,b] \subset I$. The value
\[
\ivar_p^q(x;a,b) = \sup \Bigl(\sum_{i=1}^N \bigl(\omega_q(x;t_{i-1}, t_i)\bigr)^p \Bigr)^{1/p},
\]
where the supremum is taken over all finite partitions $a = t_0 < t_1 < ... < t_N = b$ of the interval $[a,b]$, is called the $q$-integral $p$-variation of the function $x$. If $\ivar_p^q(x;a,b)<+\infty$, then we say that $x$ is a function of bounded  $q$-integral $p$-variation. The set of all such functions is denoted by $IBV_p^q([a,b])$. Later, instead of $IBV_p^q(I)$ we will simply write $IBV_p^q$.
\end{definition}

We may observe that the Definition \ref{def:ivar} given above may be rephrased a little bit. 
Let ${\mathcal I}=\{I_1, I_2, \ldots, I_N\}$ be a finite collection of nonoverlapping   subintervals of $[a,b]$ of the form $I_i=[a_i,b_i]$  and let $x\in L^q$. We set
\[
\sigma_{p,q}(x,{\mathcal I}) =  \Bigl(\sum_{i=1}^N \bigl(\omega_q(x;a_i, b_i)\bigr)^p \Bigr)^{1/p}.
\]

\begin{proposition}Let $p,q \in [1,+\infty)$ and let $x \colon [a,b] \to \mathbb R$ be a Lebesgue measurable function, where $[a,b] \subset I$. Then
\[
\ivar_p^q(x;a,b) = \sup_{{\mathcal I}} \sigma_{p,q}(x,{\mathcal I}),
\]
where the supremum is taken over all finite collections ${\mathcal I}$ of nonoverlapping subintervals of $[a,b]$.
\end{proposition}
\begin{proof}
That $\ivar_p^q(x;a,b) \leq  \sup_{{\mathcal I}} \sigma_{p,q}(x,{\mathcal I})$ is obvious because each finite partition forms the finite collection of nonoverlapping intervals. On the other hand, for each finite collection ${\mathcal I}$ of nonoverlapping subintervals $[a,b]$ we may find a corresponding partition $a = t_0 < t_1 < ... < t_N = b$ where for each $[a_i,b_i]\in{\mathcal I}$ there exists $j\in\{1, ..., N\}$ such that $a_i = t_{j-1}$ and $b_i = t_j$. Then
\[
 \sigma_{p,q}(x,{\mathcal I}) \leq \Bigl(\sum_{i=1}^N \bigl(\omega_q(x;t_{i-1}, t_i)\bigr)^p \Bigr)^{1/p},
\]
what proves that
\[
 \sup_{{\mathcal I}} \sigma_{p,q}(x,{\mathcal I}) \leq \ivar_p^q(x;a,b)
\]
and completes the proof of the Proposition.
\end{proof}

\subsection{Equivariated sets and compactness in $BV$.}

 We will now remind the definition of the {\em equivariated set} given in \cite{BG}.
 
\begin{definition}\label{eqivar:def}
A set $A\subseteq BV$ is said to be \emph{equivariated}, if for each $\varepsilon>0$ there exists a finite collection ${\mathcal I}_\varepsilon = \{I_1,I_2, \ldots,I_{N(\varepsilon)}\}$ of nonoverlapping   subintervals of $I$ of the form $I_i=[a_i,b_i]$ such that for every $x \in A$ we have
\begin{equation}\label{eq:uniform:BV}
 \var x  \leq \varepsilon+ \sigma(x,{\mathcal I}_\varepsilon),
\end{equation} 
where 
\begin{equation}
\sigma(x,{\mathcal I}_\varepsilon) \zdef  \sum_{i=1}^{N(\varepsilon)} \abs{x(b_i) - x(a_i)}.
\end{equation}
\end{definition}

This concept appeared to be very important as far as compactness criteria in $BV$ were concerned.

\begin{theorem}[\cite{BG}*{Theorem 1}]\label{th:compact:1}
A set $A \subseteq BV$ is relatively compact if and only if the following conditions are satisfied\textup:
\begin{enumerate}[label=\textup{(\roman*)}]
 \item\label{it:thm:comp1i} $A$ is bounded\textup;
 \item\label{it:thm:comp1ii} for every $x\in BV$ the set $A-x$ is equivariated.
\end{enumerate}
\end{theorem}

\begin{theorem}[\cite{BG}*{Theorem 2}]\label{th:compact:2}
A set $A \subseteq BV$ is relatively compact if and only if the following conditions are satisfied\textup:
\begin{enumerate}[label=\textup{(\roman*)}]
 \item\label{it:thm:comp2i} $A$  is bounded\textup;
 \item\label{it:thm:comp2ii} the set $A-A$ is equivariated.
\end{enumerate}
\end{theorem}

For other spaces of functions of bounded variation some partial results related to compactness were given before -- mainly related to compact embeddings. 
Some observations related to the  compact embeddings in the spaces of functions of bounded variation in the sense of Young were given  in the paper by Ciemnoczo\l{}owski and Orlicz \cite[Section~1.5]{CO}, see also \cite[Proposition 6]{BGK2}. Similar results for spaces $\Lambda BV$ were also proved in  \cite[Proposition 5]{BGK2} and \cite[ Theorem 18, Corollary 2]{BCGS}.
We should also mention here the papers concerning the compactness in the subspace $CBV$ of continuous functions of bounded variation in the sense of Jordan: the first one by Prus-Wi\'sniowski \cite{PW} and the second one -- using Bernstein polynomials approximation -- by Czudek \cite{C}.

\section{Compactness in $\Lambda BV$}

First of all we should note that the concept of equivariated set is not working directly in the $\Lambda BV$ setting. 

\begin{example}Assume $\lambda_2 > \lambda_1$. Let $x_0,x_1\colon I\to\br$ be given by
\[
x_0(t) = \begin{cases} 2 & t = 0 \cr 1 & t \in (0,1) \cr  0 & t=1,
\end{cases}
\]
\[
x_1(t) =  \begin{cases}  0 & t = 0 \cr 1 & t \in (0,1) \cr  0 & t=1,
\end{cases}
\]

Then the two-element set $A=\{x_0, x_1\}\subset \Lambda BV$ does not satisfy the appropriately modified Definition \ref{eqivar:def}. It is easy to note that if there are two intervals $[0,a]$ and $[b,1]$ belonging to the collection ${\mathcal I}$ then the greatest value of the finite sum generated by the collection ${\mathcal I}$ would be
\[
\sigma_\Lambda(x_0,{\mathcal I}) = \frac{1}{\lambda_1} + \frac{1}{\lambda_2} < \frac{2}{\lambda_1} =\var_\Lambda(x_0).
\] 
Hence we may not get arbitrarily close to the value $\var_\Lambda(x_0)$ with such cover. So the only collection that we may take is ${\mathcal I} = \{[0,1]\}$. But with this collection we have $\sigma_\Lambda(x_1,{\mathcal I}) = 0$ and we will not get arbitrarily close to 
$\var_\Lambda(x_1) = \frac{1}{\lambda_1} + \frac{1}{\lambda_2}$.
\end{example}

Therefore we need to update our definition of an {\em equivariated set}.

\begin{definition}\label{eqivar:lambda:def}A set $A\subset \Lambda BV$ is said to be \emph{$\Lambda$-equivariated}, if for each $\varepsilon>0$ there exists a finite collection ${\mathcal I}_\varepsilon = \{I_1,I_2, \ldots,I_{N(\varepsilon)}\}$ of closed subintervals (not necessarily nonoverlapping) of $I$ of the form $I_i=[a_i,b_i]$ such that for every $x \in A$, there exist such sequence of nonoverlapping intervals ${\mathcal I}_{x,\varepsilon} \subset {\mathcal I}_\varepsilon$ that
\begin{equation}\label{eq:uniform:BV}
 \var_\Lambda x  \leq \varepsilon+ \sigma_\Lambda(x,{\mathcal I}_{x,\varepsilon}).
\end{equation} 
\end{definition}

It is easy to observe that in case of $BV$ this definition is equivalent to  Definition \ref{eqivar:def}.

\begin{lemma}\label{lem:fin:L:equiv}Any finite subset $A = \{x_1, ..., x_n\}\subset \Lambda BV$ is $\Lambda$-equivariated.
\end{lemma}
\begin{proof}
Let us fix an $\varepsilon >0$. Then for each function $x_i$ there exists such a finite set ${\mathcal I}_i$ of nonoverlapping subintervals   of the interval $I$ that 
\[
 \var_\Lambda( x_i)  \leq \varepsilon+ \sigma_\Lambda(x_i,{\mathcal I}_i).
\]
Then for the set  $A$ we may take ${\mathcal I} = \bigcup_{i=1}^n {\mathcal I}_i$. Obviously ${\mathcal I}$ is finite and for any $x_i$ there exists   sequence ${\mathcal I}_i$ of intervals taken from ${\mathcal I}$ such that
\[
 \var_\Lambda( x_i)  \leq \varepsilon+ \sigma_\Lambda(x_i,{\mathcal I}_i).
\]

\end{proof}

\begin{lemma}\label{lem:lambda:equivariated}
If $A \subseteq \Lambda BV$ is a relatively compact set, then $A$ is $\Lambda$-equivariated.
\end{lemma}
\begin{proof}
Assume now, contrary to our claim, that there exits such $\varepsilon_0>0$  that for every finite collection ${\mathcal I}=\{I_1,\ldots,I_N\}$ of  closed subintervals of $I$ such function $x_{\mathcal I} \in A$ can be found that $\sigma_\Lambda(x_{\mathcal I},\tilde {\mathcal I}) + \varepsilon_0 < \var_\Lambda(x_{\mathcal I})$ for any subsequence $\tilde {\mathcal I}\subset {\mathcal  I}$ of nonoverlapping intervals.

We are going to show that in this case it is possible to construct a sequence $(x_n)_{n \in \mathbb N}$ of elements of the set $A$ such that
\[
 \var_\Lambda(x_n - x_m) \geq \frac{1}{3}\varepsilon_0 \qquad \text{for all $n\neq m$}.
\]
This clearly would contradict the fact that $A$ is relatively compact. The construction follows the induction where $x_1\in A$ is any function.

Assume now that we have chosen the functions $x_i\in A$ for $1\leq i \leq n$ in such a way that $\var_\Lambda(x_i-x_j)\geq \frac{1}{3} \varepsilon_0$ for $i\neq j$, where $i,j\in\{1,2,\ldots,n\}$. By Lemma \ref{lem:fin:L:equiv} there exists   ${\mathcal I}_n=\{I_1^n,I_2^n,\ldots,I_{N_n}^n\}$, a finite collection of closed subintervals of $I$ such that for $j=1,2,..., n$ there exists a  subsequence of nonoverlapping intervals $\tilde {\mathcal I}_j\subset {\mathcal I}_n$ such that
\[
 \var_\Lambda(x_j) -\frac{1}{3}\varepsilon_0 \leq \sigma_\Lambda(x_j,\tilde {\mathcal I}_j) \leq \var_\Lambda(x_j), \qquad j=1,\ldots,n.
\]

Moreover, by our assumption, for ${\mathcal I}_n$ there exists a function $x_{n+1}\in A$ such that
\[
\sigma_\Lambda(x_{n+1},\tilde {\mathcal I}_j) + \varepsilon_0< \var_\Lambda( x_{n+1}),
\]
for any subsequence $\tilde {\mathcal I}_j$ of nonoverlapping sub-intervals of $I$ taken from ${\mathcal I}_n$.
We will show that $\var_\Lambda(x_j-x_{n+1})\geq \frac{1}{3}\varepsilon_0$ for all $j\in\{1,\ldots,n\}$. Let us fix $j \in \{1,\ldots,n\}$. Since
\begin{align*}
\var_\Lambda(x_j-x_{n+1}) &\geq \sigma_\Lambda(x_j-x_{n+1},\tilde {\mathcal I}_j) \geq \sigma_\Lambda(x_j,\tilde {\mathcal I}_j) - \sigma_\Lambda(x_{n+1},\tilde {\mathcal I}_j) \geq\\
& \var_\Lambda(x_j) - \frac{1}{3}\varepsilon_0 - \Biggl( \var_\Lambda( x_{n+1})-\varepsilon_0\Biggr)\\
& = \frac{2}{3}\varepsilon_0  - \Biggl(  \var_\Lambda(x_{n+1}) -  \var_\Lambda( x_j)\Biggr) \geq \frac{2}{3}\varepsilon_0 -\var_\Lambda(x_{n+1}  - x_j),
\end{align*}
we get that $\var_\Lambda(x_j-x_{n+1})\geq \frac{1}{3}\varepsilon_0$.

This completes the construction of the sequence $(x_n)_{n \in \mathbb N}$ and ends the proof.
\end{proof}

\begin{theorem}\label{th:lambda:compact:1}
A set $A \subseteq \Lambda BV$ is relatively compact if and only if the following conditions are satisfied\textup:
\begin{enumerate}[label=\textup{(\roman*)}]
 \item\label{it:thm:comp1i} $A$ is bounded\textup;
 \item\label{it:thm:comp1ii} for every $x\in \Lambda BV$ the set $A-x$ is $\Lambda$-equivariated.
\end{enumerate}
\end{theorem}

\begin{proof}
If  $A \subseteq \Lambda BV$ is relatively compact, then clearly it is bounded. Moreover, given $x \in \Lambda BV$, the set $A-x$ is also relatively compact, and so, by Lemma~\ref{lem:lambda:equivariated}, we infer that $A-x$ is $\Lambda$-equivariated. This shows that   conditions~\ref{it:thm:comp1i} and~\ref{it:thm:comp1ii} are necessary.

Now assume that the set $A$ satisfies assumptions \ref{it:thm:comp1i} and \ref{it:thm:comp1ii}. Let $(x_n)_{n \in \mathbb N}$ be an arbitrary sequence of elements of  $A$. Let us observe that since $A\subset \Lambda BV$ is bounded, then by Helly's selection theorem (see \cite[Theorem 3.2]{HLW} and the general Theorem 2.49 in \cite{ABM}) there exists a subsequence $(x_{n_k})$ of $(x_n)$, which is pointwise convergent to some function $x_0\in \Lambda BV$.

Now we will show that $x_{n_k} \to x_0$ with respect to the $\Lambda BV$-norm.

For $\varepsilon>0$, by the $\Lambda$-equivariance of the set $A-x_0$, let us choose a finite collection $\mathcal I_{\varepsilon}=\{I_1,\ldots,I_{N(\varepsilon)}\}$ of closed subintervals of $I$ of the form $I_i=[a_i,b_i]$ such that
\[
 \var_\Lambda (x_{n_k}-x_0) \leq \frac{1}{4}\varepsilon + \sigma_\Lambda(x_{n_k}-x_0,\mathcal I_{\varepsilon,k}),
\]
for certain sequence ${\mathcal I}_{\varepsilon,k} = (I_1^k, ..., I_{N(\varepsilon,k)}^k)$ taken from ${\mathcal I}_{\varepsilon}$, where $I_i^k = [a_{k,i},b_{k,i}]$.
Since, $x_{n_k} \to x_0$ pointwise on $I$, we can find a number $k_0 \in \mathbb N$ such that for all $k \geq k_0$, we have $\abs{x_{n_k}(0)-x_0(0)}\leq \frac{1}{4}\varepsilon$ and
\[
 \abs{x_{n_k}(a_{i})-x_0(a_{i})}\leq \frac{\varepsilon}{4\sum_{j=1}^{N(\varepsilon)} \frac{1}{\lambda_j}} \quad \text{and} \quad \abs{x_{n_k}(b_{i})-x_0(b_{i})}\leq \frac{\varepsilon}{4\sum_{j=1}^{N(\varepsilon)} \frac{1}{\lambda_j}}, \quad i=1,\ldots,N(\varepsilon). 
\] 
Therefore
\begin{align*}
 \sigma_\Lambda(x_{n_k}-x_0,{\mathcal I}_{\varepsilon,k}) & = \sum_{i=1}^{N(\varepsilon,k)}\frac{\babs{(x_{n_k}-x_0)(b_{k,i})-(x_{n_k}-x_0)(a_{k,i})}}{\lambda_i}\\
& \leq \sum_{i=1}^{N(\varepsilon,k)}\frac{\babs{x_{n_k}(b_{k,i})-x_0(b_{k,i})}}{\lambda_i} + \sum_{i=1}^{N(\varepsilon,k)}\frac{\babs{x_{n_k}(a_{k,i})-x_0(a_{k,i})}}{\lambda_i}\\
 \leq \frac{\varepsilon}{4\sum_{i=1}^{N(\varepsilon)} \frac{1}{\lambda_i}} 2  \sum_{i=1}^{N(\varepsilon,k)} \frac{1}{\lambda_i} \leq \frac{1}{2}\varepsilon
\end{align*}
for all $k \geq k_0$. Thus
\[
 \abs{x_{n_k}(0) -x_0(0)} +\var_\Lambda (x_{n_k}-x_0) \leq \varepsilon \qquad \text{for $k \geq k_0$},
\]
which shows that $\norm{x_{n_k}-x_0}_{\Lambda BV} \to 0$ as $k \to +\infty$ and proves that the set $A$ is relatively compact in $\Lambda BV$.
\end{proof}

By a small change in the proof we get the following theorem:

\begin{theorem}\label{th:lambda:compact:2}
A set $A \subseteq \Lambda BV$ is relatively compact if and only if the following conditions are satisfied\textup:
\begin{enumerate}[label=\textup{(\roman*)}]
 \item\label{it:thm:comp2i} $A$ is bounded\textup;
 \item\label{it:thm:comp2ii}the set $A-A$ is equivariated.
\end{enumerate}
\end{theorem}
\begin{proof}The necessity part is the same as in the proof of the Theorem \ref{th:lambda:compact:1}. As far as sufficiency is concerned we will also follow the lines of Theorem \ref{th:lambda:compact:1}.
Similarly as before, let us assume that the set $A$ satisfies assumptions \ref{it:thm:comp2i} and \ref{it:thm:comp2ii}. Then for any sequence $(x_n)_{n \in \mathbb N}$  of elements of the set $A$   there exists a subsequence $(x_{n_k})$ of $(x_n)$, which is pointwise convergent to some function $x_0\in \Lambda BV$. We are going to show that $(x_{n_k})_{k\in\bn}$ is a Cauchy sequence in $\Lambda BV$, so  $\norm{x_{n_k}-x_0}_{\Lambda BV} \to 0$ as $k \to +\infty$.

Let us fix $\varepsilon>0$. Then by the equivariance of $A-A$ we may take the finite collection  $\mathcal I_{\varepsilon}=\{I_1,\ldots,I_{N(\varepsilon)}\}$ of closed subintervals of $I$ of the form $I_i=[a_i,b_i]$ such that
\[
 \var_\Lambda (x_{n_k}-x_{n_l}) \leq \frac{1}{4}\varepsilon + \sigma_\Lambda(x_{n_k}-x_{n_l},\mathcal I_{\varepsilon,k,l}),
\]
for certain sequence ${\mathcal I}_{\varepsilon,k,l} = (I_1^{k,l}, ..., I_{N(\varepsilon,k,l)}^{k,l})$ taken from ${\mathcal I}_{\varepsilon}$, where $I_i^{k,l} = [a_{k,l,i},b_{k,l,i}]$.

Since the sequence $x_{n_k}$ converges pointwise, each of the sequences $(x_{n_k}(a_i))_{k\in\bn}$ and $(x_{n_k}(b_i))_{k\in\bn}$ for $i=1,2...,N(\varepsilon)$ is a real Cauchy sequence. Hence there exists such $k_0\in\bn$ that for all $k,l\geq k_0$ there is 
\[
 \abs{x_{n_k}(a_{i})-x_{n_l}(a_{i})}\leq \frac{\varepsilon}{4\sum_{j=1}^{N(\varepsilon)} \frac{1}{\lambda_j}} \quad \text{and} \quad \abs{x_{n_k}(b_{i})-x_{n_l}(b_{i})}\leq \frac{\varepsilon}{4\sum_{j=1}^{N(\varepsilon)} \frac{1}{\lambda_j}}, \quad i=1,\ldots,N(\varepsilon). 
\] 
Then
\[
 \sigma_\Lambda(x_{n_k}-x_{n_l},\mathcal I_{\varepsilon,k,l}) \leq  \sum_{i=1}^{N(\varepsilon,k,l)}\frac{\babs{x_{n_k}(b_{k,l,i})-x_{n_l}(b_{k,l,i})}}{\lambda_i} + \sum_{i=1}^{N(\varepsilon,l,k)}\frac{\babs{x_{n_k}(a_{k,l,i})-x_{n_l}(a_{k,l,i})}}{\lambda_i} \leq \frac{\varepsilon}{4} +  \frac{\varepsilon}{4}.
\]
We may also assume that 
\[
|x_{n_k}(0)-x_{n_l}(0)|\leq  \frac{\varepsilon}{4}
\]
and then $\norm{x_{n_k}-x_{n_l}}_{\Lambda BV} \leq \varepsilon$ what proves that $(x_{n_k})_{k\in\bn}$ is a Cauchy sequence in $\Lambda BV$, which completes the proof.
\end{proof}

\section{Compactness in $\Phi BV$}

The definition of an {\em equivariated set} will be slightly different than the one in $\Lambda BV$ space, but it is going to follow a similar idea.

\begin{definition}\label{eqivar:phi:def}A set $A\subset \Phi BV$ is said to be \emph{$\phi$-equivariated}, if for each $\varepsilon>0$ there exists a finite collection ${\mathcal I}_\varepsilon = \{I_1,I_2, \ldots,I_{N(\varepsilon)}\}$ of closed subintervals of $I$ of the form $I_i=[a_i,b_i]$ such that for every $x \in A$, there exist such sequence of nonoverlapping intervals ${\mathcal I}_{x,\varepsilon} \subset {\mathcal I}_\varepsilon$ that
\begin{equation}\label{eq:uniform:BV}
 V_\phi(x)  \leq \varepsilon+ S_\phi(x,{\mathcal I}_{x,\varepsilon}).
\end{equation} 
\end{definition}

The following  lemma is the obvious consequence of Definition \ref{eqivar:phi:def} and its  proof will be omitted.
\begin{lemma}\label{lem:fin:Phi:equiv}Any finite subset $A = \{x_1, ..., x_n\}\subset \Phi BV$ is $\phi$-equivariated.
\end{lemma}

The next lemma may be proved by repeating the same steps as in the proof of Lemma \ref{lem:lambda:equivariated}, so its proof will be omitted.

\begin{lemma}\label{lem:phi:equivariated}
If $A \subseteq \Phi BV$ is a relatively compact set, then $A$ is $\phi$-equivariated.
\end{lemma}

Now, the compactness characterization in $\Phi BV$ spaces is given along the same lines as in $\Lambda BV$ spaces. 

\begin{theorem}\label{th:phi:compact:1}
A set $A \subseteq \Phi BV$ is relatively compact if and only if the following conditions are satisfied\textup:
\begin{enumerate}[label=\textup{(\roman*)}]
 \item\label{it:thm:comp:phi1i} $A$ is bounded\textup;
 \item\label{it:thm:comp:phi1ii} for every $x\in \Phi BV$ the set $A-x$ is $\phi$-equivariated.
\end{enumerate}
\end{theorem}
\begin{proof}
The necessity part is a consequence of  Lemma \ref{lem:phi:equivariated} and some general facts related to the compactness in the normed spaces (see the proof of Theorem \ref{th:lambda:compact:1} above).

The sufficiency part, similarly as before, starts with a version of Helly's selection principle (see \cite[Section 1.3]{MuOr}  and the general Theorem 2.49 in \cite{ABM}):  for any sequence $(x_n)_{n \in \mathbb N}$  of elements of the set $A$    there exists the subsequence $(x_{n_k})_{k\in\bn}$ pointwise convergent to a function $x_0\in\Phi BV$. 
Now we are going to show that $x_{n_k} \to x_0$ with respect to the $\Phi BV$-norm.

Let us fix $\varepsilon>0$. By the $\phi$-equivariance of $A-x_0$ there exists   a finite collection $\mathcal I_{\varepsilon}=\{I_1,\ldots,I_{N(\varepsilon)}\}$ of closed subintervals of $I$ of the form $I_i=[a_i,b_i]$ such that
\[
 V_\phi (x_{n_k}-x_0) \leq \frac{1}{3}\varepsilon + S_\phi(x_{n_k}-x_0,\mathcal I_{\varepsilon,k}),
\]
for certain sequence ${\mathcal I}_{\varepsilon,k} = (I_1^k, ..., I_{N(\varepsilon,k)}^k)$ taken from ${\mathcal I}_{\varepsilon}$, where $I_i^k = [a_{k,i},b_{k,i}]$. We may also assume that there exists such $\eta>0$, that 
\begin{equation}\label{eq:phi:cond:eta}
N(\varepsilon) \phi\left(\frac{2\eta}{\varepsilon / 3}\right) \leq 1.
\end{equation}
Let us assume that $k_0\in\bn$ is such that for $k\geq k_0$ the following relations hold true:
\begin{itemize}
\item $|x_{n_k}(0)-x_0(0)|\leq \frac{\varepsilon}{3}$;
\item $|x_{n_k}(a_i)-x_0(a_i)|\leq \eta$ and $|x_{n_k}(b_i)-x_0(b_i)|\leq \eta$ for all $i=1,...,N(\varepsilon)$.
\end{itemize}
Then 
\[
\sigma_\phi\left( \frac{x_{n_k}-x_0}{\varepsilon / 3},\mathcal I_{\varepsilon,k}\right) = \sum_{i=1}^{N(\varepsilon,k)} \phi\left( \frac{|x_{n_k}(b_{k,i})-x_0(b_{k,i}) - (x_{n_k}(a_{k,i})-x_0(a_{k,i}))|}{\varepsilon / 3} \right) \leq
\]
\[
 \sum_{i=1}^{N(\varepsilon,k)} \phi\left( \frac{|x_{n_k}(b_{k,i})-x_0(b_{k,i})|+| (x_{n_k}(a_{k,i})-x_0(a_{k,i}))|}{\varepsilon / 3} \right) \leq N(\varepsilon)  \phi\left(\frac{2\eta}{\varepsilon / 3}\right) \leq 1.
\]
This implies that 
\[
S_\phi( x_{n_k}-x_0,\mathcal I_{\varepsilon,k}) \leq \frac{\varepsilon}{3}
\]
That's why, for $k\geq k_0$, we have
\[
\Vert  x_{n_k}-x_0\Vert_{\Phi BV} = | x_{n_k}(0)-x_0(0)| + V_\phi( x_{n_k}-x_0) \leq \frac{\varepsilon}{3} + S_\phi( x_{n_k}-x_0,\mathcal I_{\varepsilon,k})  + \frac{\varepsilon}{3} \leq \varepsilon,
\]
which completes the proof.
\end{proof}

Similarly as in Theorem \ref{th:lambda:compact:2} we may give an alternative characterization of compactness in $\Phi BV$. The proof will be omitted since it is a modification of the proof of Theorem \ref{th:phi:compact:1} following the same ideas that were used in the proof of Theorem \ref{th:lambda:compact:2}. 
The second version of the compactness characterization is given as
\begin{theorem}\label{th:Phi:compact:2}
A set $A \subseteq \Phi BV$ is relatively compact if and only if the following conditions are satisfied\textup:
\begin{enumerate}[label=\textup{(\roman*)}]
 \item\label{it:thm:comp:phi2i} $A$ is bounded\textup;
 \item\label{it:thm:comp:phi2ii}the set $A-A$ is $\phi$-equivariated.
\end{enumerate}
\end{theorem}

\section{Compactness  in $IBV_p^q$}

From now on we are going to assume that $1\leq p < q < +\infty$. In the sequel we will refer to the well-known compactness result in $L^q(\br^n)$ spaces

\begin{theorem}[cf. Theorem 5 in \cite{HOH}]\label{th:Kolm:Riesz}Let $q\in[1,+\infty)$. A subset $A$ of $L^q(\br^n)$ is relatively compact if and only if
\begin{enumerate}
\item[(i)] $A$ is bounded;
\item[(ii)] for every $\varepsilon>0$ there is some $R>0$ so that, for every $x\in A$,
\[
\int_{|s|\geq R} |x(s)|^q\textup ds < \varepsilon^q,
\]
\item[(iii)] for every $\varepsilon>0$ there is some $\delta>0$ so that, for every $x\in A$ and $h\in\br^n$ with $|h|<\delta$,
\[
\int_{\br^n}|x(s+h) - x(s)|^q\textup ds < \varepsilon^q.
\]
\end{enumerate} 
\end{theorem}

Since we are focused on spaces $L^q(J)$, where $J=[a,b]\subset\br$ is a compact interval, we will refer to a simple corollary from the above Theorem.
\begin{corollary}\label{cor:Kolm:Riesz}A subset $A$ of $L^q(J)$ is relatively compact if and only if
\begin{enumerate}
\item[(i)] $A$ is bounded;
\item[(ii)] for every $\varepsilon>0$ there is some $\delta>0$ such that, for every $x\in A$ and $h>0$ with $|h|<\delta$,
\[
\int_{a}^{b-h}|x(s+h) - x(s)|^q\textup ds < \varepsilon^q.
\]
\end{enumerate} 
\end{corollary}

First we are going to prove an important property of $IBV_p^q$, which makes this space quite different from the spaces with variation defined {\em pointwise} (like $BV$, $\Lambda BV$ or $\Phi BV$).

\begin{lemma}\label{lem:absc:ibv}For all $x\in IBV_p^q$ we have an  absolute continuity of the $q$-integral $p$-variation with respect to the measure of the interval i.e.
\[
\forall_{\varepsilon>0} \exists_{\delta >0} \forall_{[a,b]\subset I}  b-a\leq \delta \Rightarrow \ivar_p^q(x,[a,b])\leq \varepsilon.
\]
\end{lemma} 
\begin{proof}Assume, contrary to our claim, that for some $\varepsilon_0>0$ there exists the sequence of such intervals $[\tilde a_n,\tilde b_n]\subset I$ that $\tilde b_n-\tilde a_n \leq \frac{1}{n}$ and $\ivar_p^q(x,[\tilde a_n,\tilde b_n])\geq \varepsilon_0^{1/p}$. Due to the compactness of $I$  we may assume that $\tilde a_n\to a$ and $\tilde b_n\to a$ for some $a\in I$.

We may also assume that both sequences are monotone: if we take
\[
a_n = \inf_{k\geq n}\tilde a_n
\]
and
\[
b_n = \sup_{k\geq n}\tilde b_n,
\]
then $(a_n)$ is nondecreasing and $(b_n)$ nonincreasing, $a_n\to a$, $b_n\to a$ and $[\tilde a_n,\tilde b_n]\subset[a_n,b_n]$ what implies $\ivar_p^q(x,[ a_n, b_n])\geq \ivar_p^q(x,[\tilde a_n,\tilde b_n]) \geq \varepsilon_0^{1/p}$.

Let $\delta > 0$ be such that 
\[
|A|\leq \delta \Rightarrow \left(\int_A |x(t)|^q\textup dt\right)^{1/q} \leq \frac{1}{2}\left(\frac{\varepsilon_0}{8}\right)^{1/p}.
\]
Then if $b_n-a_n\leq \delta$, then $\omega_q(x,[a_n,b_n]) \leq  \left(\frac{\varepsilon_0}{8}\right)^{1/p}$ (cf. Equation \eqref{eq:Lqmod:est} above).
Let us take such a partition $\pi_n = (t_n^0, t_n^1, ..., t_n^{k_n})$ of the interval $[a_n,b_n]$ that $\sigma_{p,q}(x,\pi_n)^p\geq \ivar_p^q(x,[ a_n, b_n])^p - \varepsilon_0/4$.  For a fixed partition the point $a$ belongs to one or two intervals forming the partition. In any case there exists such $j\in\{1,...,k_n-1\}$, that $a\in[t_n^{j-1},t_n^{j+1}]$ and $a\not\in[t_n^{i-1},t_n^{i}]$
for $i\not\in\{j,j+1\}$. This implies that 
\[
\sigma_{p,q}(x,\pi_n)^p - \omega_q(x,[t_n^{j-1},t_n^{j}])^p - \omega_q(x,[t_n^{j},t_n^{j+1}])^p \geq \ivar_p^q(x,[ a_n, b_n])^p -  \frac{\varepsilon_0}{4} - \frac{\varepsilon_0}{8} - \frac{\varepsilon_0}{8}= 
\]
\[
\ivar_p^q(x,[ a_n, b_n])^p - \frac{\varepsilon_0}{2}.
\]

Let us now take such $m\geq n$ that $[a_m,b_m] \subset [t_n^{j-1},t_n^{j+1}]$. Then  $\ivar_p^q(x,[ a_m, b_m])^p \geq \varepsilon_0^p$ and
\[
\ivar_p^q(x,[ a_n, b_n])^p \geq \ivar_p^q(x,[ a_n, t_n^{j-1}])^p  + \ivar_p^q(x,[ t_n^{j-1},t_n^{j+1}])^p + \ivar_p^q(x,[ t_n^{j+1}, b_n])^p\geq
\]
\[
\sum_{i=1}^{j-1}  \omega_q(x,[t_n^{i-1},t_n^{i}])^p + \ivar_p^q(x,[ a_m, b_m])^p +  \sum_{i=j+2}^{k_n}  \omega_q(x,[t_n^{i-1},t_n^{i}])^p =
\]
\[ \sigma_{p,q}(x,\pi_n)^p - \omega_q(x,[t_n^{j-1},t_n^{j}])^p - \omega_q(x,[t_n^{j},t_n^{j+1}])^p + \ivar_p^q(x,[ a_m, b_m])^p \geq 
\]
\[
 \ivar_p^q(x,[ a_n, b_n])^p - \varepsilon_0/2 + \varepsilon_0 = \ivar_p^q(x,[ a_n, b_n])^p + \varepsilon_0/2.
\]
a contradiction. This completes the proof.
\end{proof}

\begin{definition}\label{eqivar:ibv:def}
A set $A\subseteq IBV_p^q$ is said to be \emph{$I_p^q$-equivariated}, if for each $\varepsilon>0$ there exists a finite collection ${\mathcal I}_\varepsilon = \{I_1,I_2, \ldots,I_{N(\varepsilon)}\}$ of subintervals of $I$ of the form $I_i=[a_i,b_i]$ such that for every $x \in A$, we have a finite subcollection ${\mathcal I}_{\varepsilon,x}  \subset {\mathcal I}_\varepsilon$ of such nonoverlapping subintervals of $I$ that
\begin{equation}\label{eq:uniform:BV}
 \ivar_p^q x  \leq \varepsilon+ \sigma_{p,q}(x,{\mathcal I}_{\varepsilon,x}).
\end{equation} 
\end{definition}

The next Lemma is similar to Lemma \ref{lem:fin:L:equiv} and its obvious proof will be omitted.
\begin{lemma}\label{lem:fin:ibv:equiv}Any finite subset $A = \{x_1, ..., x_n\}\subset IBV_p^q$ is $I_p^q$-equivariated.
\end{lemma}

The proof of the next lemma is a simple modification of the proof of the Lemma \ref{lem:lambda:equivariated} and will be omitted.
\begin{lemma}\label{lem:ibv:equivariated}
If $A \subseteq IBV_p^q$ is a relatively compact set, then $A$ is $I_p^q$-equivariated.
\end{lemma}

Let us now refer to the Example 5 from \cite{JG:ibv}. This example was supposed to show that the embedding $IBV_1^q \subset L^q$ is not completely continuous. The same example shows that being equivariated and bounded in $IBV_p^q$ does not suffice for compactness.

First, let us remind the simple and important proposition

\begin{proposition}[\cite{JG:ibv}*{Example 2}]\label{prop:ex:fun:ibv}
Let $x\in L^q$ be the step function given by
\begin{equation}\label{x:step}
x(t) = \begin{cases} u, & t\in[0,c), \cr w, & t\in [c,1], \end{cases}
\end{equation}

 Then $x\in IBV_p^q$ and 
\begin{equation}
\ivar_p^q(x) = |u-w|\bigl(\min(c,1-c)\bigr)^{1/q} = \omega_q(x;0,1).
\end{equation}
\end{proposition}

\begin{example}\label{ex:not:eqv}Let the sequence  $(x_n)_{n\in\bn}$ of functions belonging to $L^q$ be given by

\begin{equation}\label{def:ex:xn}
x_n(t) = \begin{cases} n^{1/q}, & t\in [0,\frac{1}{n}], \cr 0, & t\in(\frac{1}{n},1].   \end{cases}
\end{equation}
Here $\Vert x_n\Vert_{L^q} = 1$ and $\ivar_p^q(x_n;0,1) = n^{1/q}\cdot \frac{1}{n^{1/q}} = 1$. Therefore, by Proposition \ref{prop:ex:fun:ibv},    this sequence is uniformly bounded
 in $IBV_p^q$. But the sequence is not relatively compact in $L^q$, hence not relatively compact in $IBV_p^q$. 
Because
\[
\Bigl( \int_0^{1-\frac{1}{n}} |x_n(t+\frac{1}{n}) - x_n(t)|^q \textup dt \Bigr)^{1/q}  = 1 \hfil \text{ for $n\geq 2$},
\]
for $\varepsilon = \frac{1}{2}$ there is no $\delta>0$ for which the condition (ii) of Corollary \ref{cor:Kolm:Riesz} is satisfied, which contradicts compactness.

At the same time, the sequence is $I_p^q$-equivariated: let us take ${\mathcal I} =\{ [0,1]\}$. We can see that
\[
\ivar_p^q(x_n) =  \omega_q(x_n; [0, 1]).
\]
\end{example}

As we mentioned above, the space $IBV_p^q$ is not compactly embedded in $L^q$.
It appears though, that we may show that the space $IBV_p^q$ is compactly embedded in $L^1$.

\begin{lemma}\label{compact:l1}A bounded set $A \subseteq IBV_p^q$ is relatively compact as a subset of $L^1$.
\end{lemma}
\begin{proof}  We will refer to the  Kolmogorov-Riesz compactness characterization (see Corollary \ref{cor:Kolm:Riesz} above). We are going to estimate the integral
\[
\int_0^{1-h} |x(t+h)-x(t)|\textup dt
\]
for small values of $h>0$. Let $N\in\bn$ be a fixed even number  $N = 2M$ and let $h \in(0,\frac{1}{2N})$. Then
\[
\int_0^{1-h} |x(t+h)-x(t)|\textup dt \leq \sum_{k=1}^{N-1} \int_{\frac{k-1}{N}}^{\frac{k}{N}}  |x(t+h)-x(t)|\textup dt  +  \int_{\frac{N-1}{N}}^{1-h}  |x(t+h)-x(t)|\textup dt  \leq
\]
\[
\leq \sum_{k=1}^{M} \int_{\frac{2k-2}{N}}^{\frac{2k}{N}-h}  |x(t+h)-x(t)|\textup dt  +  \sum_{k=1}^{M-1} \int_{\frac{2k-1}{N}}^{\frac{2k+1}{N}-h}  |x(t+h)-x(t)|\textup dt   \leq
\]
\[
\leq \sum_{k=1}^{M} \left(\frac{2}{N}\right)^{1/q'} \left(  \int_{\frac{2k-2}{N}}^{\frac{2k}{N}-h}  |x(t+h)-x(t)|^q\textup dt\right)^{1/q} +
 \sum_{k=1}^{M-1} \left(\frac{2}{N}\right)^{1/q'} \left( \int_{\frac{2k-1}{N}}^{\frac{2k+1}{N}-h}  |x(t+h)-x(t)|^q\textup dt\right)^{1/q}  
 \]
 \[
 \leq \left(  \sum_{k=1}^{M} \left(\frac{2}{N}\right)^{p'/q'}\right)^{1/p'} \left(  \sum_{k=1}^{M}  \left(  \int_{\frac{2k-2}{N}}^{\frac{2k}{N}-h}  |x(t+h)-x(t)|^q\textup dt\right)^{p/q}\right)^{1/p} +
 \]
 \[
 \left(  \sum_{k=1}^{M-1} \left(\frac{2}{N}\right)^{p'/q'}  \right)^{1/p'} \left(  \sum_{k=1}^{M-1} \left( \int_{\frac{2k-1}{N}}^{\frac{2k+1}{N}-h}  |x(t+h)-x(t)|^q\textup dt\right)^{p/q} \right)^{1/p}
 \]
 \[
 \leq  \left(\frac{2}{N}\right)^{1/q'-1/p'}2\ivar_p^q(x),
\]
by H\"older inequality, where $1/q'+1/q = 1$ and $1/p'+1/p = 1$ for $p\geq 1$. Because $p<q$ there is $1/q'-1/p' >0$. Of course when $h\to 0$ we may assume that $N\to+\infty$, proving that the set $A$ is relatively compact in $L^1$.
\end{proof}

Before we proceed to the next example we should refer to some properties of $L^q$-modulus of continuity. These properties are probably {\em a folclore} (especially the second one), but I could not find any reference stating these facts. Hence short proofs follow.

\begin{proposition}\label{prop:omq:1}Let $x\in L^q(a,b)$ and $\eta \in(0,b-a)$ be fixed. Then
\[
\omega_q(x;a,b) \leq \omega_q(x;a+\eta,b) + 2\xi
\]
where $\xi = \sup_{c\in[a,b-\eta]} \left( \int_c^{c+\eta} |x(t)|^q\textup dt \right)^{1/q}$.
\end{proposition}
\begin{proof}We have
\[
\omega_q(x;a,b) = \sup_{h\in(0,b-a)}  \left( \int_a^{b-h} |x(t+h) - x(t)|^q\textup dt \right)^{1/q}.
\]
For $h\in [b-a-\eta, b-a)$ there is $0<b-h-a\leq \eta$ so
\[
 \left( \int_a^{b-h} |x(t+h) - x(t)|^q\textup dt \right)^{1/q} \leq  \left( \int_a^{b-h} |x(t+h)|^q\textup dt \right)^{1/q} +  \left( \int_a^{b-h} |x(t)|^q\textup dt \right)^{1/q} =  
 \]
 \[
 \left( \int_{a+h}^{b} |x(t)|^q\textup dt \right)^{1/q} +  \left( \int_a^{b-h} |x(t)|^q\textup dt \right)^{1/q} \leq 2\xi \leq \omega_q(x;a+\eta,b) + 2\xi
\]
On the other hand for $h\in(0,\eta)$ there is  
\[
 \left( \int_a^{b-h} |x(t+h) - x(t)|^q\textup dt \right)^{1/q} \leq  \left( \int_a^{a+\eta} |x(t+h) - x(t)|^q\textup dt \right)^{1/q} +  \left( \int_{a+\eta}^{b-h} |x(t+h) - x(t)|^q\textup dt \right)^{1/q} 
 \]
 \[
 \leq 2\xi + \omega_q(x;a+\eta,b),
 \]
 what completes the proof.
\end{proof}

\begin{proposition}\label{prop:omq:2}Let $x,y\in L^q(a,b)$. Then
\[
\omega_q(x+y;a,b) \leq \omega_q(x;a,b) + \omega_q(y;a,b).
\]
\end{proposition}
\begin{proof}This result is a direct consequence of the triangle inequality in the space $L^q(a,b-h)$:
\[
 \left( \int_a^{b-h} |(x+y)(t+h) - (x+y)(t)|^q\textup dt \right)^{1/q} = 
 \]
 \[
  \left( \int_a^{b-h} |x(t+h) - x(t) + y(t+h) - y(t)|^q\textup dt \right)^{1/q} \leq 
  \]
  \[
   \left( \int_a^{b-h} |x(t+h) - x(t)|^q\textup dt \right)^{1/q} +  \left( \int_a^{b-h} |y(t+h) - y(t)|^q\textup dt \right)^{1/q}.
\]
Taking the suprema completes the proof.
\end{proof}

\begin{example}Let $A\subset  IBV_1^q$ be the set as in   Example \ref{ex:not:eqv}, i.e. $A = \{ x_n : n\in\bn \}$, where $x_n\in L^q$ is given by \eqref{def:ex:xn}. We have seen that $A$ is $I_1^q$-equivariated, but we are going to show that for all $x\in  IBV_1^q$ the set $A-x$ is $I_1^q$-equivariated as well. 

Let us fix $x\in IBV_1^q$ and $\varepsilon>0$. We may  select such a family ${\mathcal I}$ of nonoverlapping, closed subintervals of $I$ that
\[
\ivar_1^q(x) \leq \sigma_{1,q}(x,{\mathcal I}) + \frac{\varepsilon}{3}.
\]
We may assume that one of the intervals from ${\mathcal I}$ is a neighbourhood of $0$, so $[0,\delta] \in{\mathcal I}$ for certain $\delta>0$ (if there is no interval being the neighbourhood of $0$ we may add it to collection ${\mathcal I}$ increasing the value of $ \sigma_{1,q}(x,{\mathcal I})$). Let us denote
\[
{\mathcal I} = \{ [0,\delta], I_1,..., I_N\},
\]
where $N\in\bn$. 

There exists such $\eta>0$ that for any interval $[a,b]\subset I$, $b-a\leq\eta$ there is $\int_a^b |x(t)|^q\textup dt \leq \left(\frac{\varepsilon}{6}\right)^q$. We may assume that $\eta<\delta$.

Let $m\in\bn$ be such that $\frac{1}{m}<\eta$. The finite set $\{ x_n-x : n=1,...,2m\}$ is $I_1^q$-equivariated, so it is enough to show that the set $\{ x_n - x : n = 2m+1,2m+2,...\}$ is $I_1^q$-equivariated as well.

Let us now  take the collection $\tilde{\mathcal I}$ with the interval $[0,\delta]$ replaced with two intervals $\{[0,\eta],[\eta,\delta]\}$, i.e.
\[
\tilde{\mathcal I} = \{ [0,\eta],[\eta,\delta], I_1,..., I_N\}.
\]
We are going to show that 
\[
\ivar_1^q(x_n-x) \leq \sigma_{1,q}(x_n-x,\tilde{\mathcal I}) + \varepsilon,
\]
for any $n\geq 2m$.

 We can see that 
\begin{equation}\label{ex3:eq1}
 \ivar_1^q(x_n - x) \leq   \ivar_1^q(x_n) + \ivar_1^q(x) =  \ivar_1^q(x_n,[0,\eta]) +  \ivar_1^q(x),
\end{equation}
 because $x_n(t)=0$ for $t\in(\frac{1}{n},1]$. Moreover $\frac{1}{n}<\frac{\eta}{2}$,  hence (cf. Proposition \ref{prop:ex:fun:ibv}) this implies 
\begin{equation}\label{ex3:eq2}
 \ivar_1^q(x_n) =  n^{1/q}\min\left\{\frac{1}{n^{1/q}}, (1-\frac{1}{n})^{1/q}\right\} = % n^{1/q}\min\left\{\frac{1}{n^{1/q}}, (\eta-\frac{1}{n})^{1/q}\right\} = 
 \omega_q(x_n;0,\eta). 
\end{equation}
 So, for $n$ big enough, there is
\begin{equation}\label{ex3:eq3}
 \ivar_1^q(x_n - x) \leq  \ivar_1^q(x_n) +  \ivar_1^q(x) \leq  \omega_q(x_n;0,\eta) +  \ivar_1^q(x) \leq 
\end{equation}
 \[
 \omega_q(x_n;0,\eta)  +  \sigma_{1,q}(x, {\mathcal I}) + \frac{\varepsilon}{3} =   \omega_q(x_n;0,\eta) +  \omega_q(x,[0,\delta]) + \sum_{i=1}^N \omega_q(x,I_i) + \frac{\varepsilon}{3} 
 \]
 By Proposition \ref{prop:omq:2} there is
 \[
  \omega_q(x_n;0,\eta) \leq  \omega_q(x_n-x;0,\eta) +  \omega_q(x;0,\eta) \leq \omega_q(x_n-x;0,\eta) + 2\left( \int_0^\eta |x(t)|^q\textup dt \right)^{1/q} \leq \omega_q(x_n-x;0,\eta) + \frac{\varepsilon}{3}.
 \]
 On the other hand, by Proposition \ref{prop:omq:1} there is
\begin{equation}\label{ex3:eq4}
 \omega_q(x;0,\delta) \leq \omega_q(x;\eta,\delta) + 2\frac{\varepsilon}{6}.
\end{equation}
We can see that for $t\in(\eta,1]$ there is $|x(t)| = |x_n(t) - x(t)|$ and $\omega_q(x_n-x,I_i) = \omega_q(x,I_i)$ for $i=1,...,N$. Then \eqref{ex3:eq1}-\eqref{ex3:eq4} imply that
 \[
  \ivar_1^q(x_n - x) \leq \omega_q(x_n-x;0,\eta) + \frac{\varepsilon}{3} + \omega_q(x_n-x;\eta,\delta) + \frac{\varepsilon}{3} +  \sum_{i=1}^N \omega_q(x_n-x,I_i)+ \frac{\varepsilon}{3} = \sigma_{1,q}(x_n-x,\tilde{\mathcal I}) + \varepsilon,
 \]
 which shows that the sequence $(x_n-x)$ is $I_1^q$-equivariated.
\end{example}

In the proof of the next theorem we will refer to the Vitali Convergence Theorem and a series of other well-known facts.

\begin{theorem}[Vitali Convergence Theorem  (see~\cite{Lojasiewicz}*{Theorem~6.2.12})]\label{th:Vitali}
Let $y,y_n \in L^1$, $n\in \bn$ and $\sup_{n\in \bn} \int_A |y_n(t)|\textup dt \to 0$ as $|A|\to 0$ for measurable $A\subset I$. If $y_n(t)\to y(t)$ for a.e. $t\in I$, then 
\[
\int_I y_n(t)\textup dt \to \int_I y(t)\textup dt.
\]
\end{theorem}

First Lemma is probably {\em a folclore}, but I was not able to find such result in the literature -- the closest probably was the Exercise  10(g) at the end of Chapter 6 of Rudin's book \cite{Rudin}.
\begin{lemma}\label{th:l1comp:eq:i}If the set $A\subset L^q$ is compact then it is uniformly $L^q$-integrable, i.e. for any $\varepsilon >0$ there exists such $\delta >0$ that 
  for any measurable set $E\subset I$, $|E|\leq \delta$ and any $x\in A$ there is $\int_E |x(t)|^q \leq \varepsilon^q$.
\end{lemma}
\begin{proof}Let us fix $\varepsilon>0$ and a finite set of points $x_1,..,x_N\in A$ such that
\[
A \subset \bigcup_{i=1}^N B_{L^q}(x_i,\varepsilon/2).
\]
Since the set $\{x_1,...,x_N\}$ is finite, the functions $|x_i|^q$ are uniformly integrable, so there exists such $\delta>0$ that for $|E|\leq \delta$
\[
\int_E |x_i(t)|^q\textup dt < \left(\frac{\varepsilon}{2}\right)^q, \hskip1cm i = 1,2,..., N.
\]
Then for any $x\in A$ we find such $i\in\{1,..., N\}$ that $\Vert x-x_i\Vert_{L^q} < \varepsilon/2$ and
\[
\left(\int_E |x(t)|^q\textup dt\right)^{1/q} \leq \left(\int_E |x_i(t)|^q\textup dt\right)^{1/q} + \left(\int_E |x(t)-x_i(t)|^q\textup dt\right)^{1/q} \leq \varepsilon
\]
for any $E\subset I$, $|E|\leq \delta$. This proves that the set $A$ is $L^q$-uniformly integrable.
\end{proof}

The following  is a well-known fact from functional analysis.

\begin{theorem}[see Section 8.3 and Theorem 7 in Section 10.2 \cite{lax}]\label{th:wk:comp} The space $L^q$, $q\in(1,+\infty)$ is reflexive, so for each bounded sequence $(x_n)_{n\in\bn}\subset L^q$ there exists a weakly convergent subsequence.
\end{theorem}

The next theorem is well-known and usually appears as one of the first steps in the proof of the completeness of   $L^q$ spaces. 
\begin{theorem}[cf. Theorem 7.23 in \cite{axler}]\label{th:almost}A  sequence $(x_n)\subset L^1$ converging to $x_0\in L^1$ (in norm) contains a subsequence $(x_{n_k})$ such that $x_{n_k}(t) \to x_0(t)$ for almost every $t\in I$.
\end{theorem}

Now, it is time to prove the main result of this section.

\begin{theorem}\label{th:ibv:compact:1}
A set $A \subseteq IBV_p^q$ is relatively compact if and only if the following conditions are satisfied\textup:
\begin{enumerate}[label=\textup{(\roman*)}]
 \item\label{it:thm:ibv:comp1i} $A$ is bounded\textup;
  \item\label{it:thm:ibv:comp1ii} $A$ is {\em uniformly integrable} in $L^q$, i.e. for any $\varepsilon >0$ there exists such $\delta >0$ that 
  for any measurable set $E\subset I$, $|E|\leq \delta$ and any $x\in A$ there is $\int_E |x(t)|^q\textup dt \leq \varepsilon^q$;
 \item\label{it:thm:ibv:comp1iii} for every $x\in IBV_p^q$ the set $A-x$ is $I_p^q$-equivariated.
\end{enumerate}
\end{theorem}

\begin{proof}
If the set $A \subseteq IBV_p^q$ is relatively compact, then clearly it is bounded. Moreover, given $x \in IBV_p^q$, the set $A-x$ is also relatively compact, and so, by Lemma~\ref{lem:ibv:equivariated}, we infer that $A-x$ is equivariated. By Lemma \ref{th:l1comp:eq:i} each compact subset of $L^q$ is uniformly $L^q$-integrable. This shows that the conditions~\ref{it:thm:ibv:comp1i},~\ref{it:thm:ibv:comp1ii} and~\ref{it:thm:ibv:comp1iii} are necessary.

Now, we are going to show that these conditions are also sufficient. Assume $A\subset IBV_p^q$ satisfies conditions~\ref{it:thm:ibv:comp1i},~\ref{it:thm:ibv:comp1ii} and~\ref{it:thm:ibv:comp1iii}. Let $(x_n)_{n \in \mathbb N}$ be an arbitrary sequence of elements of the set $A$. Clearly it is a bounded subset of $L^q$. By Theorem \ref{th:wk:comp}  there exists a subsequence $(x_{n_k})_{k\in\bn}$ weakly convergent to some function $x_0\in L^q$. Hence the sequence $(x_{n_k})_{k\in\bn}$ converges weakly to $x_0$ in $L^1$ as well.

Since the set $A$ is bounded, by Lemma \ref{compact:l1} we can see that it is relatively compact in $L^1$, hence we may assume (taking the appropriate subsequence) that the sequence $(x_{n_k})_{k\in\bn}$ selected above converges in $L^1$ to some function $y_0\in L^1$ as well. Now the sequence 
 $(x_{n_k})_{k\in\bn}$ converges weakly in $L^1$ to $y_0$ and to $x_0$. Hence $y_0=x_0\in L^q$.

Having the sequence $(x_{n_k})_{k\in\bn}$ converging to $x_0\in L^q$ in $L^1$ norm, we may assume (taking the appropriate subsequence) that $(x_{n_k})_{k\in\bn}$ converges to $x_0$ almost everywhere (by Theorem \ref{th:almost}). This observation, assumption \ref{it:thm:ibv:comp1ii} and Vitali's convergence theorem (Theorem \ref{th:Vitali} above), with $y_n=|x_n-x_0|^q$ and $y=0$, imply that $\Vert x_{n_k}-x_0\Vert_{L^q} \to 0$

Now let us observe that $x_0\in IBV_p^q$. Let us take any partition $0=t_0<t_1 < ...< t_N= 1$ and estimate the integral
\[
\int_{t_{i-1}}^{t_i - h} |x_0(t+h) - x_0(t)|^q \textup dt.
\]
We can see that
\[
\left( \int_{t_{i-1}}^{t_i - h} |x_0(t+h) - x_0(t)|^q \textup dt \right)^{1/q} \leq
\]
\[
\left( \int_{t_{i-1}}^{t_i - h} |x_0(t+h) - x_{n_k}(t+h)|^q \textup dt \right)^{1/q} +
\]
\[
\left( \int_{t_{i-1}}^{t_i - h} |x_{n_k}(t+h) - x_{n_k}(t)|^q \textup dt \right)^{1/q} +
\]
\[
\left( \int_{t_{i-1}}^{t_i - h} |x_0(t) - x_{n_k}(t)|^q \textup dt \right)^{1/q} \leq
\]
\[
\omega_q(x_{n_k};t_{i-1}, t_i) + 2 \left( \int_{t_{i-1}}^{t_i} |x_0(t) - x_{n_k}(t)|^q \textup dt \right)^{1/q}.
\]
Since this holds for any $k\in\bn$ we may take $k$ big enough so each of the integrals satisfies
\[
\int_{t_{i-1}}^{t_i} |x_0(t) - x_{n_k}(t)|^q \textup dt \leq \frac{1}{N^{q/p}}
\]
 and eventually 
\[
\omega_q(x_0; t_{i-1},t_i)^p \leq  \left( \omega_q(x_{n_k};t_{i-1}, t_i) +  \frac{2}{N^{q/p}}\right)^p \leq C(p)\omega_q(x_{n_k};t_{i-1}, t_i)^p  + C(p)\frac{2^p}{N^{q}}
\]
for some constant $C(p)$ depending on $p$.
Then $\sum_{i=1}^N \omega_q(x_0; t_{i-1},t_i)^p \leq C(p) (\ivar_p^q (x_{n_k}))^p + \frac{2^p C(p)}{N^{q-1}}$.

Now we are going to show that $x_{n_k} \to x_0$ with respect to the $IBV_p^q$-norm.

For $\varepsilon>0$, by the equivariance of the set $A-x_0$, let us choose a finite collection $\mathcal I_{\varepsilon}=\{I_1,\ldots,I_{N(\varepsilon)}\}$ of  closed subintervals of $I$ of the form $I_i=[a_i,b_i]$ such that
\begin{equation}\label{eq:ivarpg:eqvar}
\ivar_p^q (x_{n_k}-x_0) \leq \frac{1}{3}\varepsilon + \sigma_{p,q}(x_{n_k}-x_0,\mathcal I_{\varepsilon,k}),
\end{equation}
for the appropriate subset ${\mathcal I}_{\varepsilon,k}$ of  $\mathcal I_{\varepsilon}$.
Following \eqref{eq:Lqmod:est} and because $x_{n_k} \to x_0$ in $L^q$  (as $k\to+\infty$)  we can see that
\[
\lim_{k\to +\infty} \omega_q(x_{n_k}-x_0; a_i, b_i) = 0
\]
for any $i=1,...,N(\varepsilon)$.

Hence we can find a number $k_0 \in \mathbb N$ such that for all $k \geq k_0$, we have 
\begin{equation}\label{eq:ivarpq:lqconv}
\Vert x_{n_k}-x_0\Vert_{L^q}\leq \frac{1}{3}\varepsilon
\end{equation}
and
\[
\omega_q(x_{n_k}-x_0; a_i, b_i)\leq \frac{\varepsilon}{3(N(\varepsilon))^{1/p}} \quad \text{for any} \quad i=1,\ldots,N(\varepsilon). 
\] 
Therefore
\begin{equation}\label{eq:ivarpq:sigmconv}
 \sigma_{p,q}(x_{n_k}-x_0,\mathcal I_{\varepsilon,k})  \leq \left(\sum_{i=1}^{N(\varepsilon,k)} \omega_q(x_{n_k}-x_0; a_i, b_i)^p\right)^{1/p} \leq \frac{1}{3}\varepsilon
\end{equation}
for all $k \geq k_0$. Thus by \eqref{eq:ivarpg:eqvar}, \eqref{eq:ivarpq:lqconv} and \eqref{eq:ivarpq:sigmconv} we have:
\[
\norm{x_{n_k} -x_0}_{L^q} +\ivar_p^q (x_{n_k}-x_0)  \leq \varepsilon \qquad \text{for $k \geq k_0$},
\]
which shows that $\norm{x_{n_k}-x_0}_{IBV_p^q} \to 0$ as $k \to +\infty$ and proves that the set $A$ is relatively compact in $IBV_p^q$.
\end{proof}

Similarly as before, by a small change in the proof we get the following theorem

\begin{theorem}\label{th:ibv:compact:2}
A set $A \subseteq IBV_p^q$ is relatively compact if and only if the following conditions are satisfied\textup:
\begin{enumerate}[label=\textup{(\roman*)}]
 \item\label{it:thm:ibv:comp1i} $A$ is bounded\textup;
  \item\label{it:thm:ibv:comp1ii} $A$ is {\em uniformly integrable} in $L^q$, i.e. for any $\varepsilon >0$ there exists such $\delta >0$ that 
  for any measurable set $E\subset I$, $|E|\leq \delta$ and any $x\in A$ there is $\int_E |x(t)|^q \leq \varepsilon^q$;
 \item\label{it:thm:ibv:comp1iii}  the set $A-A$ is equivariated.
\end{enumerate}
\end{theorem}
The modification of the proof of Theorem \ref{th:ibv:compact:1} goes along the same lines as it happened in case of $\Lambda BV$ and $\Phi BV$ spaces: we just have to modify the last step of the last proof and show that the subsequence  $(x_{n_k})_{k\in\bn}$ converging to $x_0$ in $L^1$ is actually a Cauchy sequence in $IBV_p^q$. The modification of the proof is obvious and will be omitted.

\end{document}